\documentclass[leqno,a4paper,oneside]{amsart}
\usepackage{amsmath,amssymb,amscd,mathrsfs}
\usepackage{graphicx,subfigure,caption,comment,multirow,color}
\usepackage{wrapfig}
\usepackage{mathptmx,courier}
\usepackage[scaled=1.0]{helvet}
\title{Symmetric diagrams for all strongly invertible knots up to 10 crossings}

\author{Christoph Lamm}

\theoremstyle{plain}
  \newtheorem{theorem}{Theorem}[section]
  
  \newtheorem{proposition}[theorem]{Proposition}

\theoremstyle{definition}
  \newtheorem{definition}[theorem]{Definition}
  \newtheorem{question}[theorem]{Question}
	\newtheorem{project}[theorem]{Project}
  \newtheorem{remark}[theorem]{Remark}
  
	\newtheorem{conjecture}[theorem]{Conjecture}

\newcommand{\myS}{\mathbb{S}}

\begin{document} 

\begin{abstract}
We present a table of symmetric diagrams for strongly invertible knots up to 10 crossings,
point out the similarity of transvergent diagrams for strongly invertible knots with symmetric
union diagrams and discuss open questions. 
\end{abstract}



\captionsetup{belowskip=10pt,aboveskip=10pt}

\def\scaling{0.8}

\definecolor{myboxcolour}{gray}{0.8}

\reversemarginpar


\maketitle

\section{Introduction} \label{sec:Introduction}

A knot $K \subset \myS^3$ is said to be \textsl{strongly invertible} if there is an orientation preserving smooth
involution $h$ of $\myS^3$ such that $h(K) = K$ and $h$ reverses the orientation on $K$.

Sakuma's article `On strongly invertible knots' \cite{Sakuma} which appeared in 1986 contains a table of symmetric diagrams of 
strongly invertible knot up to crossing number 9. Strongly invertible knots currently receive a lot of attention, see
for instance \cite{Boyle}, \cite{BoyleIssa}, \cite{LobbWatson} and \cite{Snape}. Our own motivation to extend the table 
of symmetric diagrams to knots with crossing number 10 stems from the related symmetry type of symmetric unions
of knots and a project to study strongly positive amphicheiral knots.

Strongly invertible knots can be depicted as \textsl{transvergent} diagrams (in which the axis of rotation lies in the
diagram plane) or as \textsl{intravergent} diagrams (where the axis is perpendicular to the diagram plane).
The purpose of this article is to give one transvergent diagram for each (prime) strongly invertible knot with $c(K) \le 10$.

\begin{definition}
For a knot $K$ we denote, as usual, by $c(K)$ the minimal crossing number of a diagram of $K$.
If $K$ is strongly invertible we denote by $c_t(K)$ the minimal crossing number of all transvergent diagrams of $K$.
\end{definition}

Symmetric equivalence can be defined for strongly invertible knots either using the conjugacy class 
of $h$ in the symmetry group of $K$ or by symmetric diagram moves (see \cite{LobbWatson}, Section 2).

Interestingly, especially when compared to the case of symmetric unions \cite{EisermannLamm2007}, for hyperbolic 
knots there are always one or two equivalent classes (two equivalent classes occurring for knots with period 2). 
Whereas Sakuma lists diagrams showing both classes simultaneously (for knots having two of them), 
we contend with only one diagram for each knot in the main part of this article. 
However, in Appendix B we give a list of diagrams for knots with two equivalence classes.

We focus on the crossing number of the symmetric diagrams and try to find minimal (transvergent) diagrams.
In many cases symmetric diagrams exist with $c(K)$ crossings, so that $c_t(K)=c(K)$. These symmetric diagrams
with minimal crossing number are occasionally also contained in knot tables (for instance in Rolfsen's table
the knots $9_{39}$, $9_{40}$, $9_{41}$ are shown with a vertical symmetry axis). 
For our list we cannot guarantee that the crossing number of a diagram is minimal if it exceeds $c(K)$ because 
we did not exhaustively check all symmetric diagrams.

\section{Motivation} \label{sec:Motivation}

\begin{figure}[hbtp]
\centering
\includegraphics[scale=0.85]{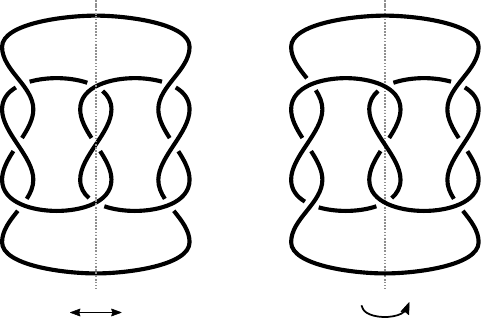}
\caption{The knot $9_{46}$. Left: as a symmetric union (mirror symmetry), right: as a strongly invertible diagram (rotational symmetry with
inversion of orientation).}
\label{9_46_st_inv_symm_union}
\end{figure}

Symmetric unions of knots were introduced by Kinoshita and Terasaka in 1957 \cite{KinoshitaTerasaka}. 
They consist of the connected sum of a knot and its mirror image with reversed orientation, with crossings inserted on the symmetry axis.
The building blocks of transvergent diagrams of strongly invertible knots are very similar: 
A connected sum of a knot diagram and a rotational copy, with crossings inserted on the symmetry axis.
Figure \ref{9_46_st_inv_symm_union} shows an example for this similarity.
In a future study we will compare methods for symmetric unions and strongly invertible knots and check if the transfer of results is possible.

\begin{figure}[hbtp]
\centering
\includegraphics[scale=1.05]{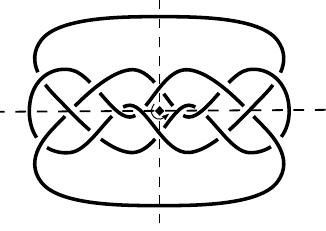}
\caption{A strongly positive amphicheiral knot diagram with two axes of symmetry (showing the  knot 12a1019)}
\label{12a1019}
\end{figure}

Another aim is the construction of strongly positive amphicheiral knots as symmetric unions with strongly invertible partial knots. 
The diagrams used for this construction possess two symmetry axes and a rotation around an axis perpendicular to the diagram plane 
yields the strongly positive amphicheiral property, see Figure \ref{12a1019}.

\enlargethispage{1cm}

A knot is called strongly positive amphicheiral if it has a diagram which is mapped 
to its mirror image by a rotation of $\pi$, preserving the orientation.
Obviously, composite knots consisting of a knot and its mirror image are strongly positive amphicheiral.
Prime strongly positive amphicheiral knots are rare, however: Until recently only three prime strongly 
positive amphicheiral knots with 12 or fewer crossings were known: $10_{99}$, $10_{123}$ and 12a427. 
In \cite{Lamm2021} we found additional examples, so that currently the list consists of $10_{99}$, $10_{123}$, 
12a427, 12a1019, 12a1105, 12a1202, and 12n706. A future study will extend this list to strongly
positive amphicheiral knots with crossing numbers up to 16 (remark: see \cite{Lamm2023}).

\section{Template notation and generation of the table} \label{sec:Notation}

All two-bridge knots are strongly invertible and because they are also 2-periodic they have two equivalent
classes of strong invertibility, see \cite{Sakuma}. Sakuma shows how to find representatives for each knot and class.
Therefore we do not need to list diagrams of two-bridge knots in detail and we concentrate on the remaining knots
up to crossing number 10 (which all have bridge number 3). For two-bridge knots we will discuss the question
whether symmetric diagrams exist with $c(K)$ crossings, however.

Our tabulating approach uses templates, in a similar way as in \cite{Lamm2021}. As half of a diagram 
contains already the information for the whole diagram, a template only shows the upper half and integer 
markings on the axis for the twists inserted on the axis. The axis is placed horizontally.

\begin{figure}[hbtp]
\centering
\includegraphics[scale=0.75]{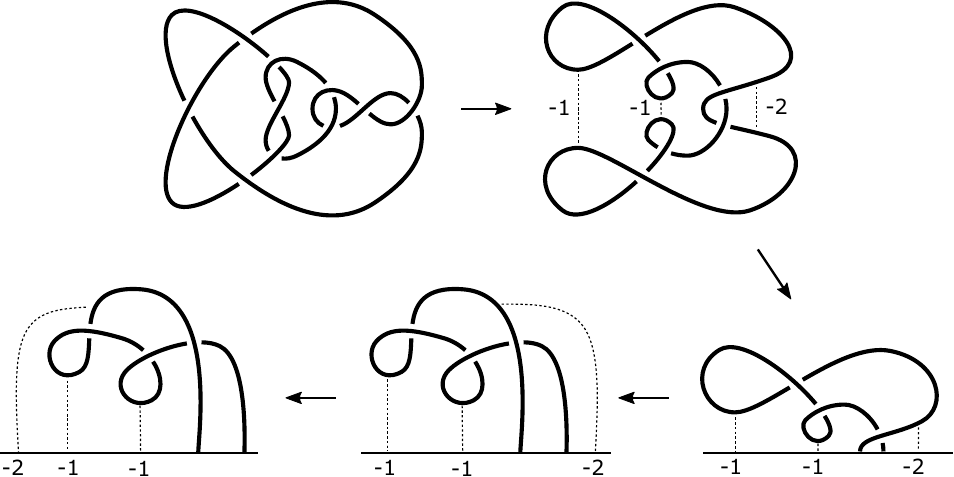}
\caption{Constructing a template representation from a transvergent diagram (the example shows $B_1(-2,-1,-1) = 10_{76}$)}
\label{diagram_to_template}
\end{figure}

Figure \ref{diagram_to_template} shows how a template diagram is constructed from a transvergent diagram.

\begin{figure}[hbtp]
\centering
\includegraphics[scale=0.8]{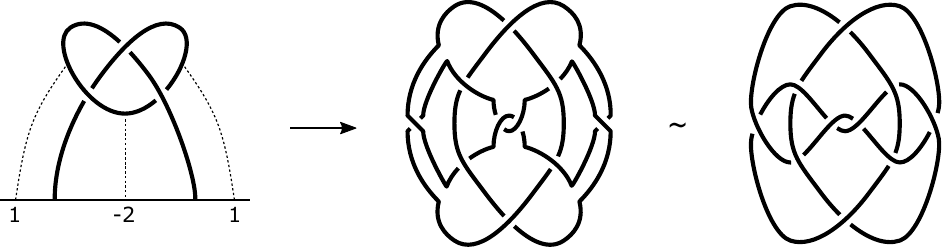}
\caption{Converting a template into a knot diagram (the example shows $C_1(1,-2,1) = 10_{122}$)}
\label{template_explained}
\end{figure}

In the other direction we proceed as in Figure \ref{template_explained}:
If a template diagram is given we need to reconstruct the other half and insert the crossings on the axis.

The list of all strongly invertible 3-bridge knots in Appendix A was compiled in the following way:
We started with Sakuma's table and converted his diagrams into templates. Variation of the twist attributes
on the axis then yielded many 10 crossing knots (using Knotscape \cite{Knotscape}). Other sources, as Rolfsen's table and 
examples in \cite{BoyleIssa} and \cite{LobbWatson} were used in the same way. The remaining 8 cases were manually 
transformed into symmetric diagrams with the help of the software KLO \cite{KLO}, see the Acknowledgments.
From the collection of all examples we chose one diagram with smallest crossing number for our table.
As we remarked, it is expected that transvergent diagrams exist with smaller crossing number for some of them.

\section{Two-bridge knots}

\begin{figure}[hbtp]
\centering
\includegraphics[scale=0.75]{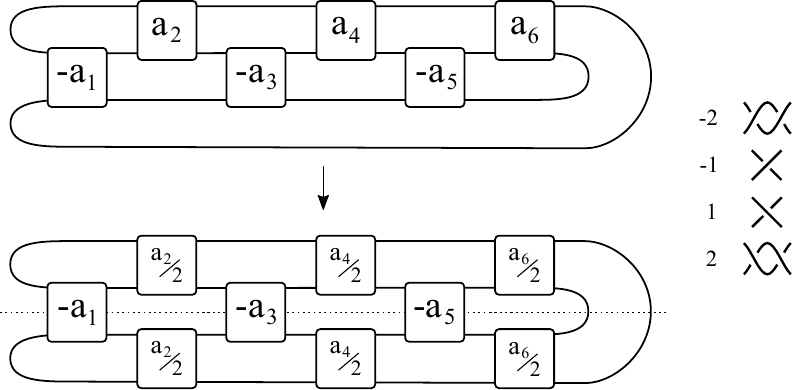}
\caption{Horizontal symmetry for alternating two-bridge knot diagrams}
\label{twobridge_def}
\end{figure}

\begin{wrapfigure}{R}{0.32\textwidth}
\centering
\includegraphics[width=0.27\textwidth]{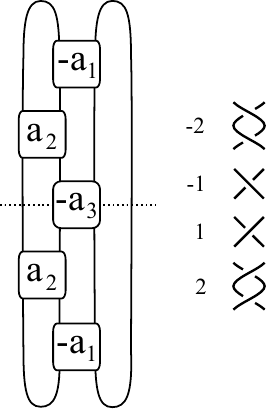}
\caption{Vertical symmetry for alternating two-bridge knot diagrams}
\label{twobridge_def_v}
\end{wrapfigure}

\begin{figure}[hbtp]
\centering
\includegraphics[scale=0.6]{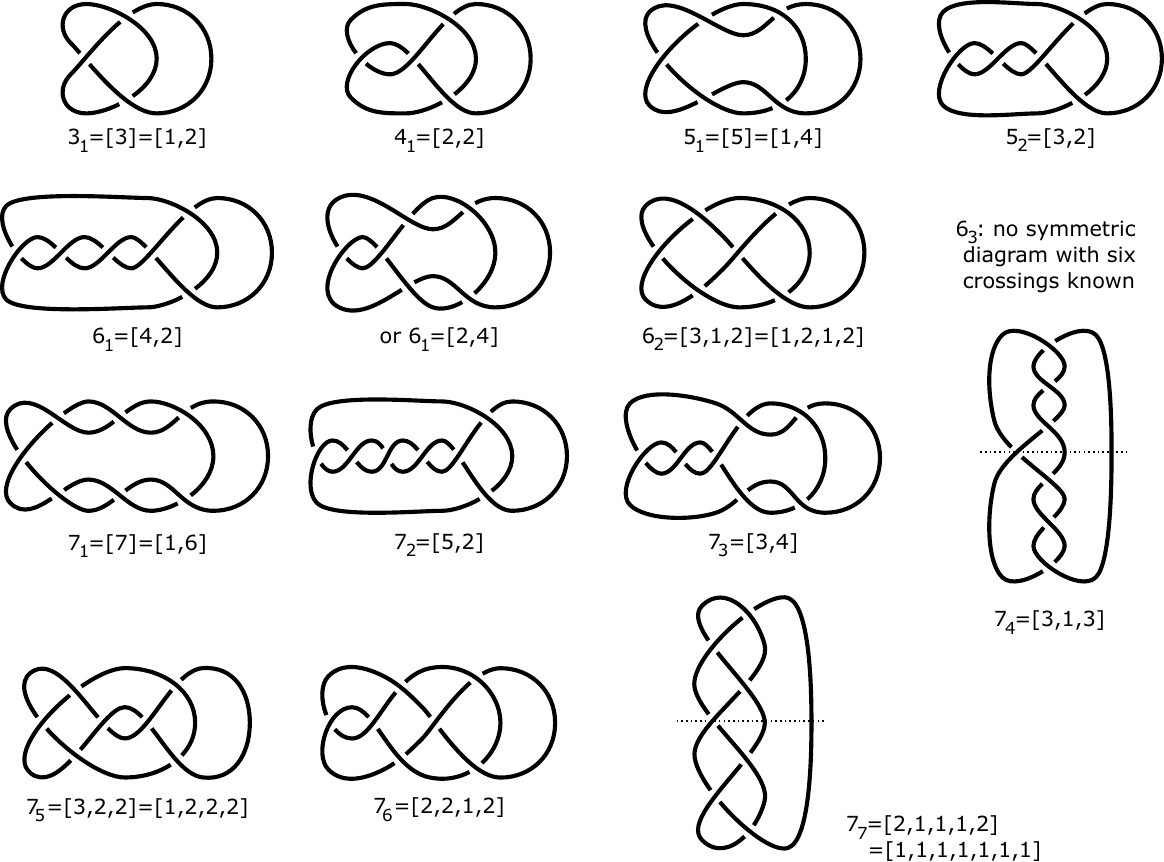}
\caption{Two-bridge knots with $c \le 7$: Symmetric diagrams with minimal crossing number}
\label{twobridge}
\end{figure}

For two-bridge knots we use the Conway notation $C(a_1, \ldots, a_n)$ (and a shorter version $[a_1, \ldots, a_n]$ in Figure \ref{twobridge}).
As usual, the plat diagram's closing patterns are different for even and odd $n$, see Figures \ref{twobridge_def} and \ref{twobridge_def_v}.

In these figures we give examples for $n=6$ (horizontal placement and axis) and $n=5$ (vertical case, but still with horizontal axis).

Sakuma's representatives for the two equivalence classes use `even' continued fractions. Each two-bridge knot
can be written as $C(a_1, \ldots, a_n)$ with even numbers $a_i$. In this case $n$ is also even and half of $n$ is the knot's genus.
If all $a_i$ are even then the flyping modification shown in Figure \ref{twobridge_def} is possible. However, to proceed as in the figure,
only ever second half-twist number needs be even.

In order to find diagrams with $c_t(K)=c(K)$ we therefore try to transform the Conway notation of each 2-bridge knot so that it realizes 
the minimal crossing number and has even entries for ever second $a_i$. 
For example, for 2-bridge torus knots we may use $T(2,2n+1) = C(2n+1) = C(1,2n)$. This condition is the first part of the following proposition. 
The second part uses a symmetric diagram which is placed vertically and is illustrated in Figure \ref{twobridge_def_v}.

\begin{proposition}
A two-bridge knot $K$ has a transvergent diagram with $c(K)$ crossings, if
\begin{itemize}
	\item $K$ can be written as $C(a_1, \ldots, a_n)$ with even $n$ and all $a_i > 0$ where the twist-numbers
	      $a_2$, $a_4$, \ldots are even, or
	\item $K$ can be written as $C(a_1, \ldots, a_n)$ with odd $n$ and all $a_i > 0$ and symmetric twist-numbers
	      $a_1 = a_n$, $a_2 = a_{n-1}$, \ldots. In this case $a_{\frac{n+1}{2}}$ is necessarily odd.
\end{itemize}
\label{twobridgeprop}
\end{proposition}

\begin{conjecture}
A two-bridge knot $K$ has a transvergent diagram with $c(K)$ crossings if and only if the above conditions are fulfilled.
\end{conjecture}\label{two-bridge-conjecture}

For crossing numbers less than 8 we illustrate this in Figure \ref{twobridge}. Note, that for the knot $6_3 = C(2,1,1,2)$ neither 
of the two conditions of Proposition \ref{twobridgeprop} is realized so that we conjecture that no transvergent diagram with 
6 crossings exists for this knot.\footnote{This should be easy to show, enumerating all alternating diagrams with 6 crossings.
Because of the Tait conjecture for alternating knots the proof of Conjecture \ref{two-bridge-conjecture} shouldn't be too difficult either.
An extended project would then try to characterize 2-bridge knots with $c_t(K) = c(K)+1$.}

There are five 2-bridge knots with 8 crossings which do not satisfy one of the two conditions:
$8_7 = C(4,1,1,2)$, $8_8 = C(2,3,1,2)$, $8_9 = C(3,1,1,3)$, $8_{13} = C(3,1,1,1,2)$, $8_{14} = C(2,2,1,1,2)$.

Remark: The article \cite{Nanasawa} by Jundai Nanasawa, published in April 2023, addresses this topic, following 
the release of the first version of the present work in October 2022.

\section{Three-bridge knots}

The result of the generation of transvergent diagrams as described in Section \ref{sec:Notation} is shown in Table \ref{tab:knotList}.
The template names require some explanation: We use $B$, $C$, $D$ and $E$ to denote template families.
Family $B$ collects curves which can be simplified by Reidemeister 1 moves. (The even simpler Family $A$ is used for the
standard representations of two-bridge knots and is described in Appendix B.) Curves in Family $C$ contain a trefoil or, as in case $C_4$,
a related diagram with a trefoil shadow. Family $D$ features the figure eight knot as half-diagram and Family $E$ the shadow of 
the torus knot $T(2,5)$. In all three cases of Family $E$ the curve is a trefoil because the crossings are not alternating.

\small
\noindent
\begin{table}[htbp]
\begin{tabular}{|l|llr|l|llr|l|llr|l|llr|}
\hline
8 a& $8_{5}$   & $B_1$ & 8  & 10 a& $10_{46}$          & $B_4$    & 10 & 10 a&  $10_{77}$   & $B_1$    & 11 & 10 n& $10_{132}$   & $C_{3a}$ & 11 \\
   & $8_{10}$  & $B_1$ & 9  &     & $10_{47}$          & $B_4$    & 11 &     &  $10_{78}$   & $C_{3a}$ & 10 &     & $10_{133}$   & $D_1$    & 11 \\
   & $8_{15}$  & $C_1$ & 8  &     & $10_{48}$          & $B_3$    & 11 &     &  $10_{89}$   & $D_1$    & 10 &     & $10_{134}$   & $B_3$    & 11 \\
   & $8_{16}$  & $C_1$ & 8  &     & $10_{49}$          & $B_3$    & 12 &     &  $10_{96}$   & $E_1$    & 13 &     & $10_{135}$   & $B_3$    & 12 \\
   & $8_{18}$  & $C_1$ & 8  &     & $10_{50}$          & $B_1$    & 10 &     &  $10_{97}$   & $E_2$    & 13 &     & $10_{136}$   & $C_1$    & 10 \\
   &           &       &    &     & $10_{51}$          & $B_1$    & 11 &     &  $10_{99}$   & $D_1$    & 10 &     & $10_{137}$   & $B_2$    & 10 \\
8 n& $8_{19}$  & $C_1$ & 8  &     & $10_{52}$          & $B_1$    & 11 &     &  $10_{100}$  & $C_{3a}$ & 10 &     & $10_{138}$   & $D_1$    & 10 \\
   & $8_{20}$  & $B_2$ & 8  &     & $10_{53}$          & $B_1$    & 12 &     &  $10_{101}$  & $C_{3a}$ & 11 &     & $10_{139}$   & $C_4$    & 10 \\
   & $8_{21}$  & $C_1$ & 8  &     & $10_{54}$          & $B_4$    & 11 &     &  $10_{103}$  & $C_1$    & 10 &     & $10_{140}$   & $B_5$    & 10 \\
   &           &       &    &     & $10_{55}$          & $B_4$    & 12 &     &  $10_{104}$  & $C_{2a}$ & 10 &     & $10_{141}$   & $C_{3b}$ & 10 \\
9 a& $9_{16}$  & $B_1$ & 9  &     & $10_{56}$          & $B_4$    & 12 &     &  $10_{105}$  & $C_{2a}$ & 11 &     & $10_{142}$   & $C_1$    & 10 \\
   & $9_{22}$  & $B_1$ & 10 &     & $10_{57}$          & $B_4$    & 13 &     &  $10_{108}$  & $C_1$    & 10 &     & $10_{143}$   & $B_5$    & 11 \\
   & $9_{24}$  & $B_1$ & 10 &     & $10_{58}$          & $B_1$    & 10 &     &  $10_{111}$  & $C_1$    & 10 &     & $10_{144}$   & $C_1$    & 10 \\
   & $9_{25}$  & $B_1$ & 10 &     & $10_{59}$          & $B_1$    & 11 &     &  $10_{112}$  & $C_{3a}$ & 10 &     & $10_{145}$   & $C_{2a}$ & 11 \\
   & $9_{28}$  & $C_6$ & 10 &     & $10_{60}$          & $D_1$    & 10 &     &  $10_{113}$  & $C_{3a}$ & 11 &     & $10_{146}$   & $C_{3a}$ & 11 \\
   & $9_{29}$  & $E_2$ & 12 &     & $10_{61}$          & $C_1$    & 10 &     &  $10_{114}$  & $C_1$    & 10 &     & $10_{152}$   & $E_3$    & 11 \\
   & $9_{30}$  & $B_1$ & 11 &     & $10_{62}$          & $B_5$    & 11 &     &  $10_{116}$  & $D_1$    & 10 &     & $10_{154}$   & $E_3$    & 11 \\
   & $9_{34}$  & $C_1$ & 9  &     & $10_{63}$          & $C_1$    & 10 &     &  $10_{120}$  & $C_1$    & 10 &     & $10_{155}$   & $C_{3b}$ & 10 \\
   & $9_{35}$  & $C_1$ & 9  &     & $10_{64}$          & $C_{2a}$ & 10 &     &  $10_{121}$  & $C_1$    & 10 &     & $10_{156}$   & $C_{2a}$ & 11 \\
   & $9_{36}$  & $B_1$ & 9  &     & $10_{65}$          & $B_5$    & 12 &     &  $10_{122}$  & $C_1$    & 10 &     & $10_{157}$   & $C_{2b}$ & 10 \\
   & $9_{37}$  & $C_1$ & 9  &     & $10_{66}$          & $C_{3a}$ & 10 &     &  $10_{123}$  & $D_1$    & 10 &     & $10_{158}$   & $C_1$    & 10 \\
   & $9_{38}$  & $E_1$ & 12 &     & $10_{68}$          & $C_{2a}$ & 11 &     &              &          &    &     & $10_{159}$   & $C_{3b}$ & 10 \\
   & $9_{39}$  & $C_1$ & 9  &     & $10_{69}$          & $C_{3a}$ & 11 & 10 n&  $10_{124}$  & $B_6$    & 10 &     & $10_{160}$   & $C_1$    & 10 \\
   & $9_{40}$  & $C_1$ & 9  &     & $10_{70}$          & $B_1$    & 10 &     &  $10_{125}$  & $B_3$    & 11 &     & $10_{161}$   & $C_4$    & 10 \\
   & $9_{41}$  & $C_1$ & 9  &     & $10_{71}$          & $B_1$    & 11 &     &  $10_{126}$  & $B_3$    & 10 &     & $10_{162}$   & $C_1$    & 10 \\
   &           &       &    &     & $10_{72}$          & $B_1$    & 11 &     &  $10_{127}$  & $B_3$    & 11 &     & $10_{163}$   & $C_1$    & 10 \\
9 n& $9_{42}$  & $C_1$ & 9  &     & $10_{73}$          & $B_1$    & 12 &     &  $10_{128}$  & $C_1$    & 10 &     & $10_{164}$   & $C_1$    & 10 \\
   & $9_{43}$  & $B_2$ & 9  &     & $10_{74}$          & $C_5$    & 10 &     &  $10_{129}$  & $B_2$    & 10 &     & $10_{165}$   & $C_1$    & 10 \\
   & $9_{44}$  & $B_2$ & 9  &     & $10_{75}$          & $C_6$    & 10 &     &  $10_{130}$  & $B_2$    & 10 &     &              &          & \\
   & $9_{45}$  & $C_1$ & 10 &     & $10_{76}$          & $C_{2a}$ & 10 &     &  $10_{131}$  & $B_2$    & 11 &     &              &          & \\
   & $9_{46}$  & $C_1$ & 9  &     &                    &          &    &     &              &          & &        &              &          & \\
   & $9_{47}$  & $C_1$ & 9  &     &                    &          &    &     &              &          & &        &              &          & \\
   & $9_{48}$  & $C_1$ & 9  &     &                    &          &    &     &              &          & &        &              &          & \\
   & $9_{49}$  & $C_1$ & 9  &     &                    &          &    &     &              &          & &        &              &          & \\
\hline 
\end{tabular}
\caption{A list of all 118 prime strongly invertible knots with bridge number 3. The second column indicates the 
template where the diagram information is found and the third column contains the crossing number of the diagram.}
\label{tab:knotList}
\end{table}

\normalsize

As remarked before, we chose one diagram for each knot for this list. For instance, the knot $10_{76}$ is assigned to template $C_{2a}$
with a symmetric diagram with 10 crossings, $10_{76}  = C_{2a}(0,0,-2)$. We could also have assigned it to template $B_1$, because, as
shown in Figure \ref{diagram_to_template}, $10_{76} = B_1(-2,-1,-1)$ in a 10 crossing diagram.

\section{A comparison of symmetries}

\begin{figure}[hbtp]
\centering
\includegraphics[scale=0.42]{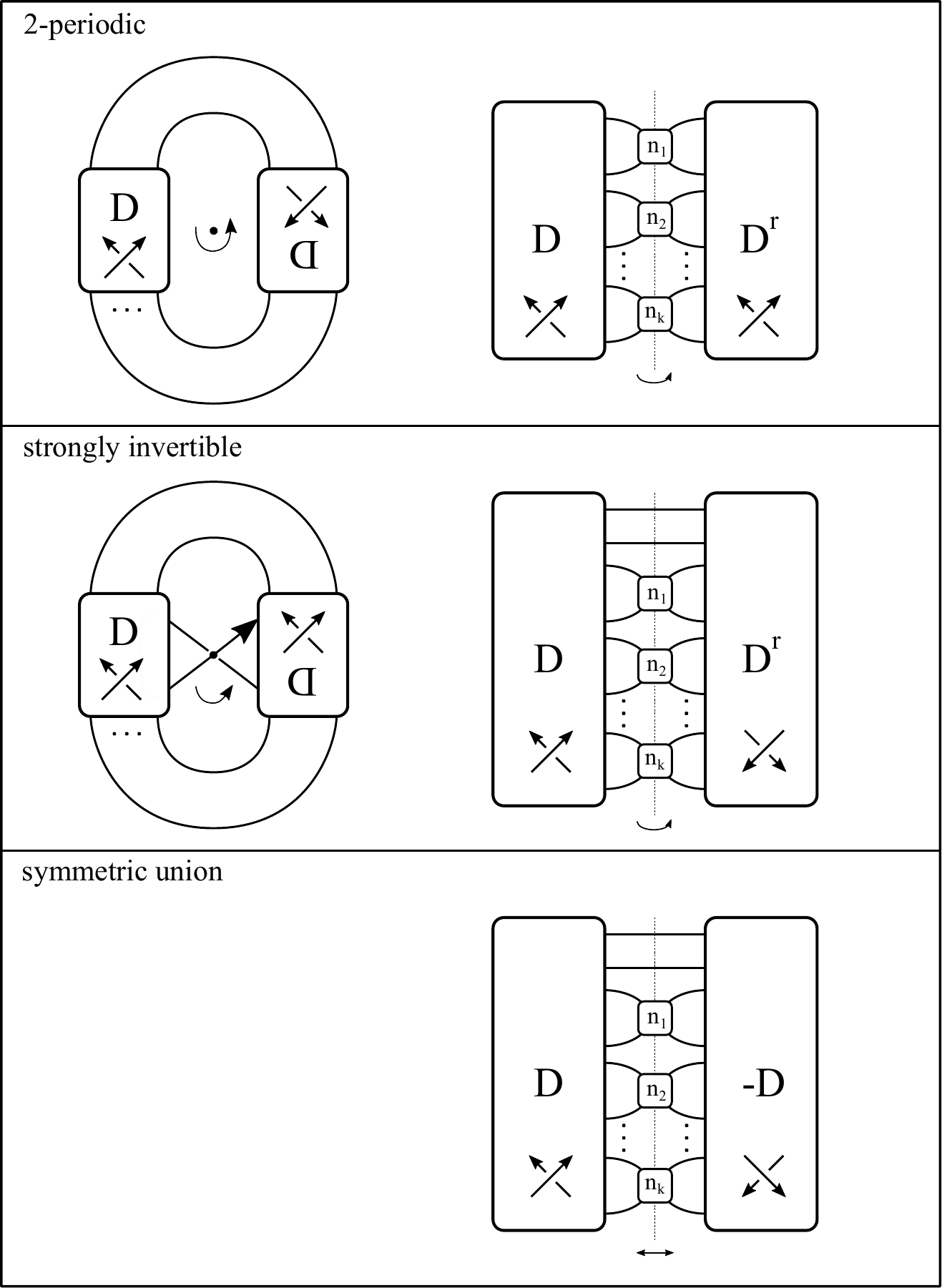}
\caption{Three symmetries: Intravergent diagrams are shown on the left and transvergent diagrams on the right. Symmetric unions are symmetric to a plane of reflection and all diagrams are transvergent.}
\label{symmetries_periodic_invertible_union}
\end{figure}

\begin{figure}[hbtp]
\centering
\includegraphics[scale=0.42]{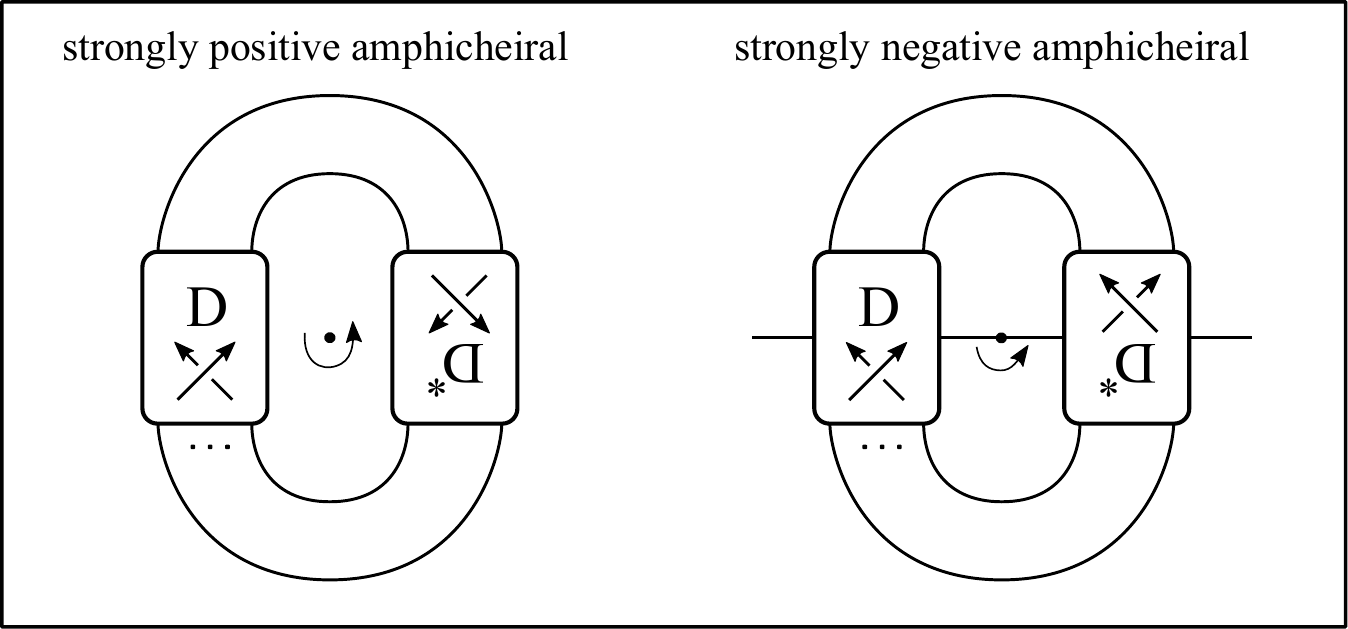}
\caption{Two amphicheiral symmetries: Strongly positive and negative amphicheiral knots are point symmetric and all diagrams are intravergent.}
\label{symmetries_amphicheiral}
\end{figure}

In this section we compare several involutive symmetries of knots, see Figures \ref{symmetries_periodic_invertible_union}
and \ref{symmetries_amphicheiral}. For a diagram $D$ we denote by $D^r$ the rotated diagram, by $-D$ the mirrored diagram (at a
plane perpendicular to the diagram plane) with reversed orientation, and by $D^*$ the diagram mirrored in the diagram plane.

\newpage

\section{Comparison of ribbons spanning symmetric unions and strongly invertible knots}

\begin{figure}[hbtp]
\centering
\includegraphics[scale=0.8]{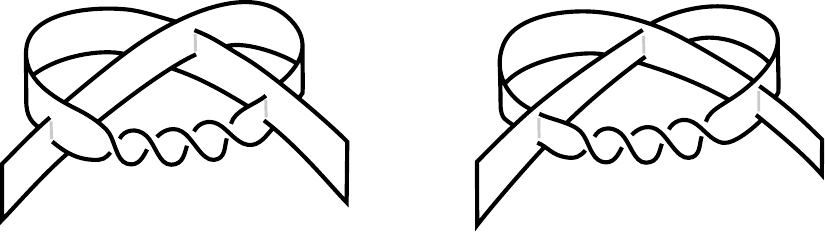}
\caption{A comparison of the singularities in a band (a) left: spanning a strongly invertible knot, $10_{142}=C_1(0,4,0)$, with
clasp singularities and (b) right: a symmetric union, $10_{140}$, with ribbon singularities.}
\label{band_diagram}
\end{figure}

As a first topic in the comparison of methods for symmetric unions and strongly invertible knots we look at the symmetric 
ribbon which spans a symmetric union \cite{Lamm2021}. Because one half in a corresponding diagram for a strongly invertible 
knot differs from the symmetric union just by switching crossings we obtain a similar ribbon but with clasp singularities, 
see Figure \ref{band_diagram}. The clasp number of a knot $K$ is the minimum number of clasp singularities among all disks 
spanning $K$. The article \cite{KadokamiKawamura} contains a table of the clasp numbers of knots up to crossing 
number 10 (some of the values are not yet known).

Example: For $10_{142}$, shown in Figure \ref{band_diagram} on the left, we find in \cite{KadokamiKawamura} that the clasp number is 3 
and the diagram contains 3 clasp singularities. For another knot, $10_{139}$, the clasp number is 4 and it is listed in our table with 
template $C_4$ which generates 4 clasp singularities in the ribbon (because there are 4 crossings in the template).
We conclude that this knot cannot be generated with templates $B_1$, $B_2$ or $C_3$.

Therefore, transferring the use of the ribbon from the symmetric union case to the strongly invertible case may give some
interesting results and, as mentioned, we propose to do that in more detail in a future article.

\section{Plausibility check for alternating diagrams and related examples}
We distinguish alternating and non-alternating templates -- defined in the same way as for knot diagrams.
Only alternating templates can generate transvergent diagrams for alternating knots with $c(K)$ crossings.
On the other hand, non-alternating templates may generate transvergent diagrams for alternating knots, 
which then have necessarily more than $c(K)$ crossings.

\subsection{Alternating diagrams}
Table \ref{tab:knotList} and Appendix A contain 8 alternating and 10 non-alternating templates. For the 8 alternating templates
we list in Table \ref{tab:alternating_templates} the alternating knots with transvergent diagrams having $c(K)$ crossings. 
For each template there is a specific condition on the twist numbers: the diagram generated from the template is alternating 
if and only if this condition is satisfied. This is a useful plausibility check for the computational results.

\small
\begin{table}[htbp]
\begin{tabular}{|l|l|l|}
\hline
$B_1$    & $t_1 \le 0$, $t_2 \le 0$, $t_3 \le 0$              & $8_5$, $9_{16}$, $9_{36}$, $10_{50}$, $10_{58}$, $10_{70}$ \\
$B_4$    & $t_1 \le 0$, $t_2 \le 0$, $t_3 \le 0$              & $10_{46}$ \\
$C_1$    & $t_1 \ge 0$, $t_2 \le 0$, $t_3 \ge 0$              & all 18 alternating knots listed with template $C_1$ \\
$C_{2a}$ & $t_1 \ge 0$, $t_2 \le 0$, $t_3 \le 0$              & $10_{64}$, $10_{76}$, $10_{104}$ \\
$C_{3a}$ & $t_1 \ge 0$, $t_2 \le 0$, $t_3 \ge 0$, $t_4 \ge 0$ & $10_{66}$, $10_{78}$, $10_{100}$, $10_{112}$ \\
$C_5$    & $t_1 \le 0$, $t_2 \le 0$                           & $10_{74}$ \\
$C_6$    & $t_1 \ge 0$, $t_2 \ge 0$                           & $10_{75}$ \\
$D_1$    & $t_1 \ge 0$, $t_2 \ge 0$, $t_3 \le 0$, $t_4 \le 0$ & $10_{60}$, $10_{89}$, $10_{99}$, $10_{116}$, $10_{123}$ \\
\hline 
\end{tabular}
\caption{The alternating templates, the conditions on the twist numbers for alternating diagrams and the respective alternating knots}
\label{tab:alternating_templates}
\end{table}

\normalsize

\begin{remark}
There are alternating knots for which Table \ref{tab:knotList} contains non-alternating diagrams and these obviously do not satisfy 
the conditions in Table \ref{tab:alternating_templates}. For instance $10_{72} = B_1(-1,-1,3)$ is an alternating knot and listed with 
template $B_1$ as a diagram with 11 crossings. The twist numbers $t_1 = -1$ and $t_2 = -1$ satisfy the condition $t_1 \le 0$, $t_2 \le 0$ 
but the third twist number does not, since $t_3 > 0$. The similar cases $10_{47}$, $10_{54}$ and $10_{55}$ are illustrated in the next
section in Figure \ref{reducible_diagrams}.
\end{remark}

\subsection{Alternating knots with `nearly alternating diagrams'}
\begin{figure}[hbtp]
\centering
\includegraphics[scale=0.8]{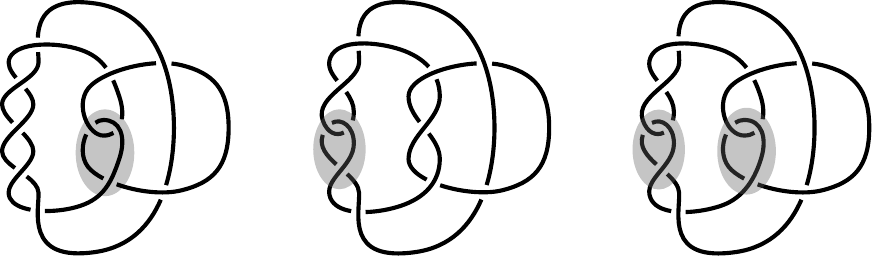}
\caption{In each case flipping the gray area(s) yields an alternating diagram: 
$10_{47} = B_4(-1,0,2)$, $10_{54} = B_4(2,0,-1)$ and $10_{55} = B_4(2,0,2)$}
\label{reducible_diagrams}
\end{figure}

Table \ref{tab:knotList} and Appendix A contain some diagrams for alternating knots with $c(K)+1$ crossings which can be transformed by
a `flip' to alternating diagrams with $c(K)$ crossings, see Figure \ref{reducible_diagrams} (and there is even the case $10_{55}$ 
with two flips, transforming a symmetric diagram with 12 crossings to a minimal alternating one). 

Since these flips are especially simple, we list the cases in Table \ref{tab:nearly_alternating}. 
In each of these cases a twist number of 2 occurs in a template where in the last section we noted the condition $t_i \le 0$, 
or a twist number of -2 where we noted the condition $t_i \ge 0$.
This corresponds to the fact that the diagrams are not alternating inside the areas marked in gray.

\small
\begin{table}[htbp]
\begin{tabular}{|l|l|l|}
\hline
$B_1$    & $8_{10}$, $9_{24}$, $9_{25}$, $9_{30}$, $10_{51}$, $10_{53}$, $10_{71}$, $10_{77}$ \\
$B_4$    & $10_{47}$, $10_{54}$, $10_{55}$ \\
$C_{2a}$ & $10_{68}$, $10_{105}$ \\
$C_{3a}$ & $10_{69}$, $10_{101}$, $10_{113}$ \\
\hline 
\end{tabular}
\caption{Symmetric diagrams in alternating templates which can be transformed by `flips' to alternating diagrams}
\label{tab:nearly_alternating}
\end{table}

\normalsize

\begin{remark}
Similar `nearly minimal' diagrams exist for non-alternating knots, for instance $10_{124}$, $10_{134}$ in the
non-alternating template $B_3$, and $10_{146}$ in the alternating template $C_{3a}$.
\end{remark}

\newpage
\section{Open questions and suggested projects}
The knot invariant $c_t(K)$ is currently only known for knots for which transvergent diagrams with $c(K)$ crossings were 
experimentally found, and in this case $c_t(K) = c(K)$. This situation can be improved in computational or theoretical ways.

\begin{project}
Generate all transvergent diagrams up to $b$ crossings and determine their knot types. 
Since for 3-bridge knots our list contains diagrams with $\le 13$ crossings, it is enough to set $b=12$ 
in order to obtain $c_t(K)$ for all (strongly invertible, prime) 3-bridge knots up to crossing number 10.
\end{project}

\begin{project}
Find a knot invariant which allows the computation of a lower bound for $c_t(K)$ for each (strongly invertible) knot $K$.
\end{project}

On the other hand, we would like to know more about the difference $c_t(K)-c(K)$:

\begin{project}
Find an upper bound for $d(n):=\max_K\{c_t(K)-c(K)|c(K) = n\}$, for strongly invertible knots $K$. 
It might be reasonable to study $d(n)$ separately for 2-bridge knots and knots with a larger bridge number.
\end{project}

\section*{Acknowledgments}
I thank Marc Kegel for his contribution of 8 symmetric diagrams in the Mathoverflow discussion 
`Diagrams for strongly invertible knots with 10 crossings' \cite{mathoverflow} which I started in May 2022.

\clearpage
\section{Appendix A}
Appendix A contains the 18 templates generating 118 diagrams of prime strongly invertible knots
with bridge number 3 and crossing number $\le 10$. The twist numbers correspond to the dotted lines
in the templates read from left to right.

\vspace{1cm}

\small
%

\hspace{-2.0cm}
\colorbox{myboxcolour}{Family B}

\noindent
\parbox[t]{4.0cm}{
\centering
\mbox{} \\
\includegraphics[scale=\scaling]{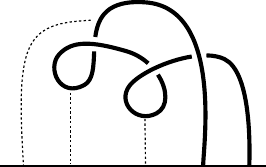} \\
template $B_1$ \\
\mbox{} \\
\begin{tabular}{rcr@{, }r@{, }r}
$8_{5\phantom{0}}$ &=&( 0& -1& -1) \\
$8_{10}$   &=&( 0& -1&  2) \\
$9_{16}$   &=&(-1& -1& -1) \\
$9_{22}$   &=&( 0& -1&  3) \\
$9_{24}$   &=&(-1& -1&  2) \\
$9_{25}$   &=&( 0& -2&  2) \\
\end{tabular}
}
\parbox[t]{4.0cm}{
\mbox{} \\
\begin{tabular}{rcr@{, }r@{, }r}
$9_{30}$   &=&( 0&  2&  3) \\
$9_{36}$   &=&( 0& -2& -1) \\
$10_{50}$  &=&( 0& -3& -1) \\
$10_{51}$  &=&( 0& -3&  2) \\
$10_{52}$  &=&( 0& -1&  4) \\
$10_{53}$  &=&( 0&  2&  4) \\
$10_{58}$  &=&( 0& -2& -2) \\
$10_{59}$  &=&( 0& -2&  3) \\
$10_{70}$  &=&(-1& -2& -1) \\
$10_{71}$  &=&(-1& -2&  2) \\
$10_{72}$  &=&(-1& -1&  3) \\
$10_{73}$  &=&(-1&  2&  3) \\
$10_{77}$  &=&(-2& -1&  2) \\
\end{tabular}
}
\noindent
\parbox[t]{3.5cm}{
\centering
\mbox{} \\
\includegraphics[scale=\scaling]{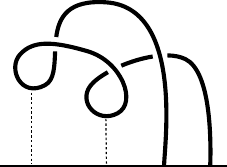} \\
template $B_2$
}
\parbox[t]{3.5cm}{
\mbox{} \\
\begin{tabular}{lcr@{, }r}
$\phantom{1}8_{20}$   &=&(-1&  1) \\
$\phantom{1}9_{43}$   &=&(-1&  2) \\
$\phantom{1}9_{44}$   &=&(-2&  1) \\
$10_{129}$ &=&(-3&  1) \\
$10_{130}$ &=&(-1&  3) \\
$10_{131}$ &=&( 2&  3) \\
$10_{137}$ &=&(-2&  2) \\
\end{tabular}
}

\vspace{-0.1cm}
\hspace{7.0cm}
\noindent
\parbox[t]{4.2cm}{
\centering
\mbox{} \\
\includegraphics[scale=\scaling]{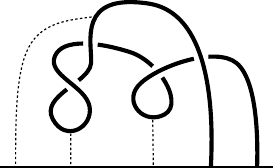} \\
template $B_3$
}
\parbox[t]{3.5cm}{
\mbox{} \\
\begin{tabular}{lcr@{, }r@{, }r}
$10_{48}$  &=&(-1&  1& -1) \\
$10_{49}$  &=&(-1&  1&  2) \\
$10_{125}$ &=&( 1&  1& -1) \\
$10_{126}$ &=&( 0&  1& -1) \\
$10_{127}$ &=&( 0&  1&  2) \\
$10_{134}$ &=&( 0& -2& -1) \\
$10_{135}$ &=&( 0& -2&  2) \\
\end{tabular}
}

\vspace{1.5cm}
\noindent
\parbox[t]{4.0cm}{
\centering
\mbox{} \\
\includegraphics[scale=\scaling]{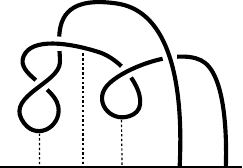} \\
template $B_4$
}
\parbox[t]{4.0cm}{
\mbox{} \\
\begin{tabular}{rcr@{, }r@{, }r}
$10_{46}$  &=&(-1&  0& -1) \\
$10_{47}$  &=&(-1&  0&  2) \\
$10_{54}$  &=&( 2&  0& -1) \\
$10_{55}$  &=&( 2&  0&  2) \\
$10_{56}$  &=&(-1&  2& -1) \\
$10_{57}$  &=&(-1&  2&  2) \\
\end{tabular}
}
\parbox[t]{3.5cm}{
\centering
\mbox{} \\
\includegraphics[scale=\scaling]{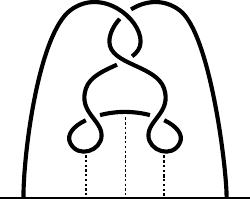} \\
template $B_5$
}
\parbox[t]{3.5cm}{
\mbox{} \\
\begin{tabular}{lcr@{, }r@{, }r}
$10_{62}$   &=&(1&  1& -1) \\
$10_{65}$   &=&(1& -2& -1) \\
$10_{140}$  &=&(1&  0& -1) \\
$10_{143}$  &=&(1& -1& -1) \\
\end{tabular}
}

\vspace{1cm}
\noindent
\parbox[t]{4.0cm}{
\centering
\mbox{} \\
\includegraphics[scale=\scaling]{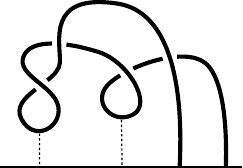} \\
template $B_6$
}
\parbox[t]{3.5cm}{
\mbox{} \\
\begin{tabular}{lcr@{, }r}
$10_{124}$ &=&( 1&  1) \\
\end{tabular}
}
\newpage

\hspace{-2.0cm}
\colorbox{myboxcolour}{Family C}

\noindent
\parbox[t]{4.0cm}{
\centering
\mbox{} \\
\includegraphics[scale=\scaling]{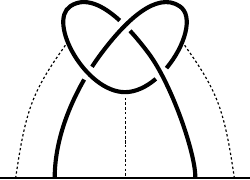} \\
template $C_1$ \\
\mbox{} \\
\begin{tabular}{lcr@{, }r@{, }r}
$\phantom{1}8_{15}$    &=&( 0&  0&  2) \\
$\phantom{1}8_{16}$    &=&( 0& -1&  1) \\
$\phantom{1}8_{18}$    &=&( 1&  0&  1) \\
$\phantom{1}8_{19}$    &=&( 0&  2&  0) \\
$\phantom{1}8_{21}$    &=&(-2&  0&  0) \\
$\phantom{1}9_{34}$    &=&( 1&  0&  2) \\
$\phantom{1}9_{35}$    &=&( 0& -3&  0) \\
$\phantom{1}9_{37}$    &=&( 0&  0&  3) \\
$\phantom{1}9_{39}$    &=&( 0& -1&  2) \\
$\phantom{1}9_{40}$    &=&( 1& -1&  1) \\
$\phantom{1}9_{41}$    &=&( 0& -2&  1) \\
$\phantom{1}9_{42}$    &=&( 0&  2&  1) \\
$\phantom{1}9_{45}$    &=&( 1&  1&  2) \\
$\phantom{1}9_{46}$    &=&( 0&  3&  0) \\
\end{tabular}
}
\parbox[t]{4.2cm}{
\centering
\mbox{} \\
\begin{tabular}{lcr@{, }r@{, }r}
$\phantom{1}9_{47}$    &=&(-2&  0&  1) \\
$\phantom{1}9_{48}$    &=&(-3&  0&  0) \\
$\phantom{1}9_{49}$    &=&(-2& -1&  0) \\
$10_{61}$   &=&( 0& -4&  0) \\
$10_{63}$   &=&( 0&  0&  4) \\
$10_{103}$  &=&( 0& -1&  3) \\
$10_{108}$  &=&( 0& -3&  1) \\
$10_{111}$  &=&( 0& -2&  2) \\
$10_{114}$  &=&( 1&  0&  3) \\
$10_{120}$  &=&( 2&  0&  2) \\
$10_{121}$  &=&( 1& -1&  2) \\
$10_{122}$  &=&( 1& -2&  1) \\
$10_{128}$  &=&( 0&  2&  2) \\
$10_{136}$  &=&( 1&  2&  1) \\
$10_{142}$  &=&( 0&  4&  0) \\
$10_{144}$  &=&(-4&  0&  0) \\
$10_{158}$  &=&(-3& -1&  0) \\
$10_{160}$  &=&( 0&  3&  1) \\
$10_{162}$  &=&(-2& -2&  0) \\
$10_{163}$  &=&(-3&  0&  1) \\
$10_{164}$  &=&(-2& -1&  1) \\
$10_{165}$  &=&(-2&  0&  2) \\
\end{tabular}
}
\noindent
\parbox[t]{3.5cm}{
\centering
\mbox{} \\
\includegraphics[scale=\scaling]{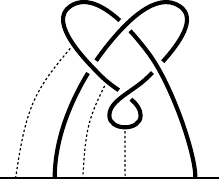} \\
template $C_{2a}$
}
\parbox[t]{3.5cm}{
\mbox{} \\
\begin{tabular}{lcr@{, }r@{, }r}
$10_{64}$  &=&( 0& -1& -1) \\
$10_{68}$  &=&( 0& -1&  2) \\
$10_{76}$  &=&( 0&  0& -2) \\
$10_{104}$ &=&( 1&  0& -1) \\
$10_{105}$ &=&( 1&  0&  2) \\
$10_{145}$ &=&( 0&  2& -1) \\
$10_{156}$ &=&( 1&  1& -1) \\
\end{tabular}
}

\vspace{-4.0cm}
\hspace{7.9cm}
\noindent
\parbox[t]{3.5cm}{
\centering
\mbox{} \\
\includegraphics[scale=\scaling]{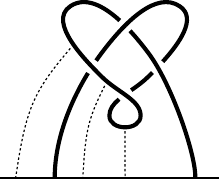} \\
template $C_{2b}$
}
\parbox[t]{3.5cm}{
\mbox{} \\
\begin{tabular}{lcr@{, }r@{, }r}
$10_{157}$  &=&( 1& 0& 1) \\
\end{tabular}
}

\vspace{2.6cm}
\hspace{-1.5cm}
\noindent
\parbox[t]{5.5cm}{
\centering
\mbox{} \\
\includegraphics[scale=\scaling]{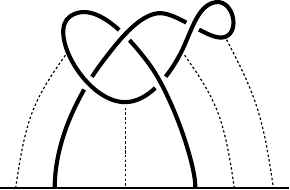} \\
template $C_{3a}$
}
\parbox[t]{4.2cm}{
\mbox{} \\
\begin{tabular}{lcr@{, }r@{, }r@{, }r}
$10_{66}$  &=&( 0&  0&  1&  1) \\
$10_{69}$  &=&( 0&  0&  1& -2) \\
$10_{78}$  &=&( 0&  0&  0&  2) \\
$10_{100}$ &=&( 0& -1&  0&  1) \\
$10_{101}$ &=&( 0& -1&  0& -2) \\
$10_{112}$ &=&( 1&  0&  0&  1) \\
$10_{113}$ &=&( 1&  0&  0& -2) \\
$10_{132}$ &=&( 0& -1& -1&  1) \\
$10_{146}$ &=&( 0&  0& -2&  1) \\
\end{tabular}
}
\noindent
\parbox[t]{5.5cm}{
\centering
\mbox{} \\
\includegraphics[scale=\scaling]{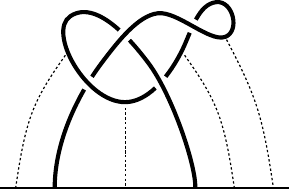} \\
template $C_{3b}$ \\
\mbox{} \\
\begin{tabular}{lcr@{, }r@{, }r@{, }r}
$10_{141}$  &=&( 0&  0& -1& -1) \\
$10_{155}$  &=&( 0& -1&  0& -1) \\
$10_{159}$  &=&( 1&  0&  0& -1) \\
\end{tabular}
}

\vspace{1.0cm}
\noindent
\parbox[b]{4.0cm}{
\centering
\mbox{} \\
\includegraphics[scale=\scaling]{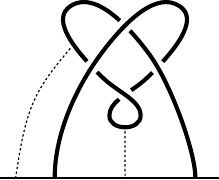} \\
template $C_4$ \\
\mbox{} \\
\begin{tabular}{rcr@{, }r}
$10_{139}$ &=&( 1& 1) \\
$10_{161}$ &=&(-1& 1) \\
\end{tabular}
}
\hspace{0.8cm}
\parbox[b]{3.8cm}{
\centering
\mbox{} \\
\includegraphics[scale=\scaling]{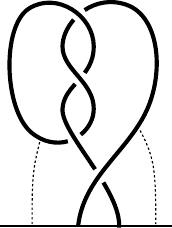} \\
template $C_5$ \\
\mbox{} \\
\begin{tabular}{rcr@{, }r}
$10_{74}$ &=&(-1&-1) \\
\end{tabular}
}
\hspace{0.8cm}
\parbox[b]{3.8cm}{
\centering
\mbox{} \\
\includegraphics[scale=\scaling]{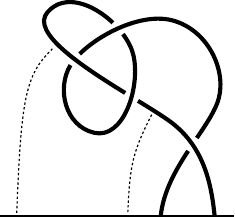} \\
template $C_6$ \\
\mbox{} \\
\begin{tabular}{rcr@{, }r}
$9_{28}$  &=&(-2& 0) \\
$10_{75}$ &=&( 1& 1) \\
\end{tabular}
}

\newpage

\hspace{-2.0cm}
\colorbox{myboxcolour}{Family D}

\noindent
\parbox[t]{4.5cm}{
\centering
\mbox{} \\
\includegraphics[scale=\scaling]{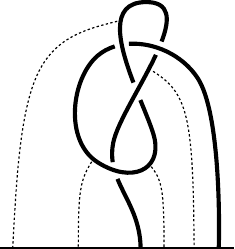} \\
template $D_1$
}
\parbox[t]{3.5cm}{
\mbox{} \\
\begin{tabular}{lcr@{, }r@{, }r@{, }r}
$10_{60}$   &=&( 2& 0& 0&  0) \\
$10_{89}$   &=&( 1& 0& 0& -1) \\
$10_{99}$   &=&( 0& 1& 0& -1) \\
$10_{116}$  &=&( 1& 1& 0&  0) \\
$10_{123}$  &=&( 1& 0&-1&  0) \\
$10_{133}$  &=&( 1& 1& 0&  1) \\
$10_{138}$  &=&( 0& 0& 2&  0) \\
\end{tabular}
}

\vspace{0.8cm}


\hspace{-2.0cm}
\colorbox{myboxcolour}{Family E}

\noindent
\parbox[t]{4.5cm}{
\centering
\mbox{} \\
\includegraphics[scale=\scaling]{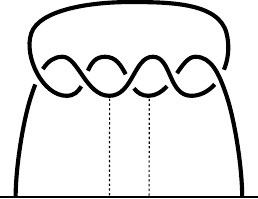} \\
template $E_1$
}
\parbox[t]{3.5cm}{
\mbox{} \\
\begin{tabular}{rcr@{, }r}
$9_{38}$  &=&(-1& 1) \\
$10_{96}$ &=&(-1& 2) \\
\end{tabular}
}

\vspace{0.8cm}
\noindent
\parbox[t]{4.5cm}{
\centering
\mbox{} \\
\includegraphics[scale=\scaling]{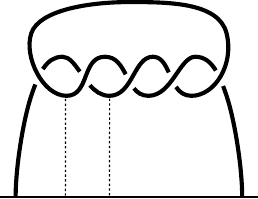} \\
template $E_2$
}
\parbox[t]{3.5cm}{
\mbox{} \\
\begin{tabular}{rcr@{, }r}
$9_{29}$  &=&( 1& 1) \\
$10_{97}$ &=&( 2& 1) \\
\end{tabular}
}

\vspace{0.8cm}
\noindent
\parbox[t]{4.5cm}{
\centering
\mbox{} \\
\includegraphics[scale=\scaling]{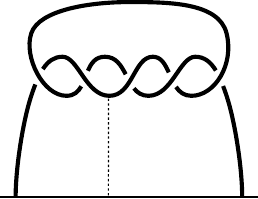} \\
template $E_3$
}
\parbox[t]{3.0cm}{
\mbox{} \\
\begin{tabular}{rcr}
$10_{152}$ &=&( 1) \\
$10_{154}$ &=&(-1) \\
\end{tabular}
}

\clearpage
\newpage
\normalsize

\section{Appendix B}
This appendix contains diagrams for strongly invertible knots with two equivalence classes of strong inversions.
As mentioned, these are knots having (cyclic or free) period 2. If possible, we use diagrams with two
perpendicular transvergent axes in the diagram plane, as shown in the following illustration:

\begin{figure}[hbtp]
\centering
\includegraphics[scale=0.8]{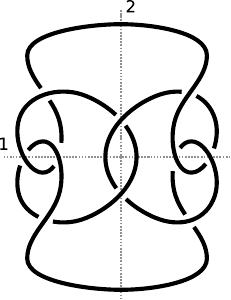}
\caption{A doubly transvergent diagram for the knot $8_{15}$.
The two axes correspond to different equivalence classes of strong inversions.}
\label{doubly_transvergent_8_15}
\end{figure}

\begin{definition}
A \textit{doubly transvergent diagram} is a transvergent knot diagram with two perpendicular axes of rotation 
in the diagram plane, each realizing a strong inversion of the knot. 
\end{definition}

In this definition, we do not require that the two strong inversions belong to different equivalence classes.
Sakuma's list of diagrams for three-bridge knots (see \cite{Sakuma}) contain the following cases as doubly transvergent
diagrams with two distinct equivalence classes: $8_5$, $8_{15}$, $8_{21}$, $9_{16}$, $9_{28}$.

The composition of the two strong inversions in a doubly transvergent diagram gives a rotation with respect to an axis 
which is perpendicular to the diagram plane, and we therefore have:

\begin{proposition} \label{prop_transvergent}
A knot with a doubly transvergent diagram is strongly invertible and has cyclic period 2.
\end{proposition}

Conversely, we show that all (prime) strongly invertible knots with cyclic period 2 up to ten crossings have doubly transvergent 
diagrams, see the next section (for two-bridge knots) and Figure \ref{template_three_bridge} (for three-bridge knots).
We then discuss, in which cases the two strong inversions belong to different equivalence classes.
We exclude torus knots and the two knots with free period 2 (these are $10_{155}$ and $10_{157}$).

For the analysis of different cases, we consider the symmetry groups.

\subsection{Symmetry groups.}
We find the following symmetry groups for prime strongly invertible knots with cyclic period 2
up to 10 crossings (excluding torus knots):

\begin{table}[htbp]
\begin{tabular}{|l|r|r|r|}
\hline
          & most frequent case & additional cases \\
\hline
2-bridge knots & $75 \times D_2$ & $(10 + 6) \times D_4$ \\
\hline
3-bridge knots & $28 \times D_2$ & $4 \times D_6$ \\
	             &                 & $1 \times D_8$ \\
\hline 
\end{tabular}
\label{tab:symmetry_groups}
\end{table}

Here, the notation $(10 + 6) \times D_4$ refers to the occurence of the symmetry group $D_4$
for amphicheiral knots (10 times) and chiral knots (6 times). We will see that their geometric
properties are quite different.

\subsection{The cases $D_2$ (these are always chiral) and $D_4$ (case of amphicheiral knots)}
The most frequent case, the dihedral group $D_2$ with four elements, occurs in the simplest situation
(prime, \textit{chiral}, strongly invertible knots with cyclic period 2). 
From the symmetry classification in \cite{BoyleRouseWilliams} (Section 14.3, type SIP) we conclude: 

\begin{proposition}
Let $K$ be a prime hyperbolic knot. If $K$ is strongly invertible, has cyclic period 2 and symmetry group $D_2$,
then it has a doubly transvergent diagram with axes belonging to different equivalence classes of strong inversions.
\end{proposition}

We have the same result for prime, \textit{amphicheiral}, strongly invertible knots with cyclic period 2 and symmetry group $D_4$
(see \cite{BoyleRouseWilliams}, Section 14.3, type SNASI): 

\begin{proposition}
Let $K$ be a prime hyperbolic knot. If $K$ is strongly invertible and amphicheiral, has cyclic period 2 and symmetry group $D_4$,
then it has a doubly transvergent diagram with axes belonging to different equivalence classes of strong inversions.
\end{proposition}

\begin{figure}[hbtp]
\centering
\includegraphics[scale=0.9]{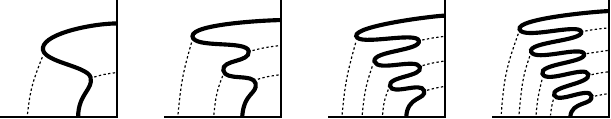}
\caption{The templates $A_1$, $A_2$, $A_3$ and $A_4$ for two-bridge knots with genus $ \le 4$.}
\label{template_two_bridge}
\end{figure}

\begin{figure}[hbtp]
\centering
\includegraphics[scale=0.9]{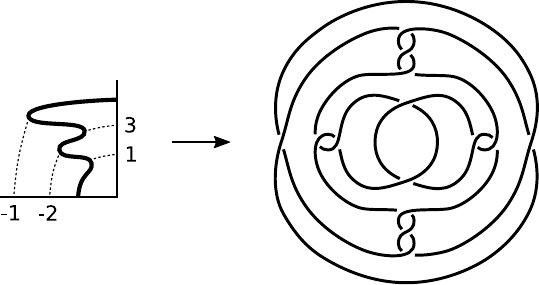}
\caption{Using template $A_2$ for a two-bridge knot with genus 2. The example $A_2(-1,-2 \mid 3,1) = C(2,6,4,2)$ is shown.}
\label{two_bridge_example}
\end{figure}

For two-bridge knots, doubly transvergent diagrams are given by Sakuma (first part of Proposition 3.6 in \cite{Sakuma}).
They can be described in template form, as shown in Figures \ref{template_two_bridge} and \ref{two_bridge_example}. 
Analogously to the template construction in Section \ref{sec:Notation} for transvergent diagrams (one axis), for doubly 
transvergent diagrams (two axes) a quarter of a diagram already contains the complete diagram information. 

Here, integer markings are necessary on both axes. The twist numbers correspond to the dotted lines in the templates.
In the template notation, as for instance in $A_2(-1,-2 \mid 3,1)$, the twist numbers are mapped to the dotted lines 
horizontally from left to right (before $\mid$), and vertically going downwards (after $\mid$).
Compare also with the construction of doubly symmetric diagrams for strongly positive amphicheiral knots in \cite{Lamm2023}.

We remark, that the 75 two-bridge (non-torus) knots with symmetry group $D_2$, together with the 10 amphicheiral two-bridge knots 
with symmetry group $D_4$, are those with fractions $p/q$ and $q^2 \not\equiv 1 (\text{mod}\;p)$.

We continue the discussion of symmetry groups and are going to see, that not in all cases doubly transvergent
diagrams with axes belonging to two different equivalence classes are possible.

\subsection{Cases with symmetry group $D_6$.}
These are the knots $9_{35}$ and $9_{40}$ with cyclic periods 2 and 3, and the knots $9_{48}$ and $10_{75}$ with cyclic period 2 and 
free period 3. We find doubly transvergent diagrams for them (see Figure \ref{template_three_bridge}), but it is currently not clear,
whether the axes belong to different equivalence classes of strong inversions. For $9_{35}$, $9_{40}$ and $9_{48}$, 
Sakuma uses transvergent diagrams with the second axis of strong inversions perpendicular to the diagram plane.
Viewing these representations from the side would give doubly transvergent diagrams, so we assume that we can find
doubly transvergent diagram with axes belonging to different equivalence classes of strong inversions in these four cases.
A general argument showing this for symmetry group $D_6$ would be desirable.

\subsection{The case with symmetry group $D_8$.}
The knot $8_{18}$ has period 4, is amphicheiral and has symmetry group $D_8$. The diagram given by Sakuma is doubly transvergent. 
The two orthogonal axes belong to the same equivalence class, however. 
An axis with different equivalence class is obtained by rotating one of the axes by $\pi/4$ (as shown for $8_{18}$ in Sakuma's diagram list. 
We assume, that this can be generalized to all (prime, strongly invertible, hyperbolic) amphicheiral knots with period 4 and symmetry group $D_8$.

\subsection{Cases with symmetry group $D_4$.}
The 6 chiral knots with symmetry group $D_4$ are the non-torus two-bridge knots with fractions $p/q$ and $q^2 \equiv 1 (\text{mod}\;p)$:

\small
\begin{tabular}{lclcl}
$7_4$    &=& $C(4,-4)$           &=& $C(3,1,3)$ \\
$7_7$	   &=& $C(2,2,-2,-2)$      &=& $C(2,1,1,1,2)$ \\
$9_{10}$ &=& $C(4,-2,2,-4)$      &=& $C(3,3,3)$ \\
$9_{17}$ &=& $C(2,2,-2,2,-2,-2)$ &=& $C(2,1,3,1,2)$ \\
$9_{23}$ &=& $C(2,-4,4,-2)$      &=& $C(2,2,1,2,2)$ \\
$9_{31}$ &=& $C(2,-2,-2,2,2,-2)$ &=& $C(2,1,1,1,1,1,2)$
\end{tabular}
\normalsize

\medskip
These knots have an additional symmetry, leading to the symmetry group of $D_4$ instead of $D_2$.
This symmetry can be seen in the even notation (middle column), or also in the positive notation (right column).
Interestingly, they coincide with the knots we earlier referred to as two-bridge knots with vertical symmetry, 
see Figure \ref{twobridge_def_v} and the second condition in Proposition \ref{twobridgeprop}.

The two axes in the standard diagram (Sakuma \cite{Sakuma}, first part of Proposition 3.6) belong to the same equivalence
class of strong inversions in these cases. Sakuma gives an additional, intravergent, diagram in the second part of his Proposition 3.6
and illustrates the two inequivalent axes in his Figure 3.2 on the right. 

\begin{figure}[hbtp]
\centering
\includegraphics[scale=0.8]{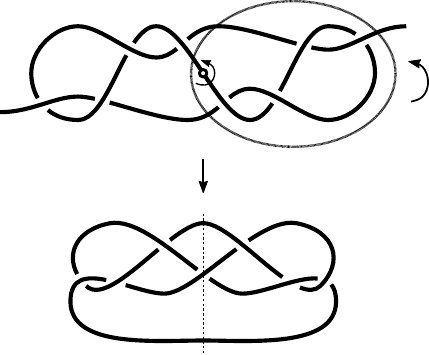}
\caption{An intravergent diagram for $7_7$, representing the second equivalence class of strong inversions (top), 
and a transformation to a transvergent diagram, with the same equivalence class (bottom).}
\label{symmetry_D_4_example_7_7}
\end{figure}

The intravergent diagram and a transformation to a transvergent diagram is shown in Figure \ref{symmetry_D_4_example_7_7} for the knot $7_7$.

\enlargethispage{0.5cm}

In Sakuma's diagram (\cite{Sakuma}, Figure 3.2 on the right), the axis for the cyclic rotation by $\pi$ is denoted by $\gamma$. The axes for the 
inequivalent strong inversions are $\alpha$ and $\beta$. Because $\alpha$ and $\gamma$ intersect in the center of the diagram, but $\beta$ does 
not intersect any of these two axes, we conclude that in these cases the situation is quite different from the previously discussed cases.

\begin{question}
Does for (non-torus) two-bridge knots with fractions $p/q$ and $q^2 \equiv 1 (\text{mod}\;p)$ exist a doubly transvergent diagram 
with the second equivalence class of strong inversions (denoted with $I_2$ by Sakuma) for both axes of inversion? 
\end{question}

We expect, that this is not possible, because the second axis $\beta$ does not intersect the other two axes.

\begin{question}
Is the geometrical situation, concerning the intersection configuration of the three axes, the same for all 
prime, hyperbolic, strongly invertible, chiral knots, with cyclic period 2 and symmetry group $D_4$?
\end{question}

\noindent
The discussion in this appendix led to a conjecture which can be found in \cite{Lamm2025}, Section 4.

\begin{figure}[hbtp]
\centering
\small
\noindent
\parbox[t]{4.0cm}{
\centering
\mbox{} \\
\includegraphics[scale=1.1]{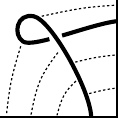} \\
template $X$ \\
\mbox{} \\
\vspace{0.85cm}
\begin{tabular}{l@{\;=\;}r@{\hspace{0cm}}r@{, }r@{, }r@{ $\mid$ }r@{, }r@{, }r}
$8_{5}$   &(& -1&  0&  0&  0&  0&  1) \\
$8_{15}$  &(&  2&  0&  0&  0&  0&  1) \\
$8_{18}$  &(& -1&  0&  0& -1&  0&  0) \\
$8_{19}$  &(& -1&  0&  0&  0&  0& -1) \\
$8_{21}$  &(&  2&  0&  0&  0&  0& -1) \\
$9_{16}$  &(& -1&  0& -1&  0&  0& -1) \\
$9_{28}$  &(&  2&  0& -1&  0&  0& -1) \\
$9_{35}$  &(&  1&  0& -1&  0&  0&  2) \\
$9_{37}$  &(&  1&  0& -1& -1&  0& -2) \\
$9_{40}$  &(&  1&  0&  0& -2&  0&  1) \\
$9_{46}$  &(&  2&  0&  0& -1&  0&  0) \\
$9_{48}$  &(&  1&  0& -1& -1&  0&  1) \\
\end{tabular}
}
\parbox[t]{4.2cm}{
\centering
\mbox{} \\
\begin{tabular}{l@{\;=\;}r@{\hspace{0cm}}r@{, }r@{, }r@{ $\mid$ }r@{, }r@{, }r}
$10_{58}$ &(& -2&  0&  0&  0&  0&  1) \\
$10_{60}$ &(&  3&  0&  0&  0&  0&  1) \\
$10_{61}$ &(& -1&  0&  0&  0&  0&  2) \\
$10_{63}$ &(&  2&  0&  0&  0&  0&  2) \\
$10_{64}$ &(& -1&  0&  1&  0&  1&  1) \\
$10_{66}$ &(&  2&  0&  1&  0&  1&  1) \\
$10_{68}$ &(&  1&  0& -1& -1&  1& -1) \\
$10_{69}$ &(&  1& -1&  1& -1& -1&  1) \\
$10_{74}$ &(&  1&  0&  1& -1&  1& -1) \\
$10_{75}$ &(&  1& -1& -1& -1& -1&  1) \\
$10_{76}$ &(& -1&  0& -1&  0&  0&  1) \\
$10_{78}$ &(&  2&  0& -1&  0&  0&  1) \\
$10_{120}$ &(& -2&  0&  0& -1&  0&  0) \\
$10_{122}$ &(& -1&  0&  0& -1&  0&  1) \\
$10_{136}$ &(&  1&  0&  0&  2&  0&  0) \\
$10_{138}$ &(& -3&  0&  0&  1&  0&  0) \\
$10_{139}$ &(& -1&  0& -1&  0& -1& -1) \\
$10_{141}$ &(&  2&  0& -1&  0& -1& -1) \\
$10_{142}$ &(& -1&  0&  0&  0&  0& -2) \\
$10_{144}$ &(&  2&  0&  0&  0&  0& -2) \\
$10_{145}$ &(&  1&  0&  1& -1&  0&  1) \\
$10_{146}$ &(&  1&  0& -2& -1&  0& -1) \\
\end{tabular}
}
\caption{Template $X$ generates all prime strongly invertible three-bridge knots 
having cyclic period 2, up to 10 crossings. This list includes the knot $8_{18}$ with
period 4 and the torus knot $8_{19}$.}
\label{template_three_bridge}
\end{figure}

\noindent
A complete list of doubly transvergent diagrams for prime three-bridge knots up 
to 10 crossings is shown in Figure \ref{template_three_bridge}.

\bigskip

\normalsize
We assume that in Figure \ref{template_three_bridge} all knots except $8_{18}$ and $8_{19}$ have two 
equivalence classes of strong inversions in the diagram described by template X (using the $x$- and $y$-axes). 
To verify this for the remaining cases $9_{35}$, $9_{40}$, $9_{48}$ and $10_{75}$, one could compute the 
$\eta$-polynomials corresponding to the two possible choices of axes.
\clearpage
\newpage


\vspace{1cm}
\noindent
Christoph Lamm \\ \noindent
R\"{u}ckertstr. 3, 65187 Wiesbaden \\ \noindent
Germany \\ \noindent
e-mail: christoph.lamm@web.de

\end{document}